\documentclass[12pt]{amsart}
\usepackage {amssymb,amsfonts,amsmath,amsopn,amstext,amscd,latexsym,xy}

\input xy
\xyoption{all}

\theoremstyle{remark}

\newtheorem{para}{\bf}[section]

\theoremstyle{definition}

\theoremstyle{plain}

\newtheorem{thm}[para]{Theorem}

\newtheorem{cor}[para]{Corollary}
\newtheorem{prop}[para]{Proposition}

\newcommand{\al}{{\alpha}}

\newcommand{\vpi}{\varpi}

\newcommand{\cC}{{\mathcal C}}

\newcommand{\cM}{{\mathcal M}}

\newcommand{\cR}{{\mathcal R}}

\newcommand{\bX}{{\bf X}}

\newcommand{\bbC}{{\mathbb C}}

\newcommand{\bbF}{{\mathbb F}}

\newcommand{\bbQ}{{\mathbb Q}}

\newcommand{\bbX}{{\mathbb X}}

\newcommand{\bbZ}{{\mathbb Z}}

\newcommand{\fronr}{{\hat{\fro}^{nr}}}
\newcommand{\hFnr}{{\hat{F}^{nr}}}

\newcommand{\Qlb}{{\overline{\mathbb{Q}_l}}}

\newcommand{\froxB}{{{\mathfrak o}_B^{\times}}}

\renewcommand{\frm}{{\mathfrak m}}

\newcommand{\fro}{{\mathfrak o}}
\newcommand{\frp}{{\mathfrak p}}

\newcommand{\End}{{\rm End}}

\newcommand{\Aut}{{\rm Aut}}

\newcommand{\hra}{\hookrightarrow}

\newcommand{\lra}{\longrightarrow}

\newcommand{\ra}{\rightarrow}

\newcommand{\Sp}{{\rm Sp}}
\newcommand{\Spec}{{\rm Spec}}
\newcommand{\Spf}{{\rm Spf}}
\newcommand{\sub}{\subset}

\begin{document}

\begin{center}\Large

{\bf Geometrically connected components of Lubin-Tate deformation
spaces with level structures}\\

\medskip \normalsize


\bigskip

{\bf Matthias Strauch}\\

\bigskip
\smallskip

{\it Department of Pure Mathematics and Mathematical Statistics\\
Centre for Mathematical Sciences, University of Cambridge\\
Wilberforce Road, Cambridge, CB3 0WB, United Kingdom\\
e-mail: M.Strauch@dpmms.cam.ac.uk}\\

\end{center}

\bigskip
\bigskip

{\small {\bf Abstract.} We determine the geometrically connected
components of the generic fibre of the deformation space $\cM_m$
which parameterizes deformations of a one-dimensional formal
$\fro$-module equipped with Drinfeld level-$m$-structures. It is
shown that the geometrically connected components are defined over
a Lubin-Tate extension of the base field, and the action of the
covering group $GL_n(\fro/\vpi^m))$ on the components is given by
the determinant. This furnishes a description of the action of
this group on the \'etale cohomology of the spaces $\cM_m^{rig}$
in degree
zero.\\

\bigskip

\normalsize

\tableofcontents

\section{Introduction}

Let $F$ be a local non-Archimedean field, and denote by $\fro$ its
ring of integers. Let $\bbX$ be a one-dimensional formal
$\fro$-module of $F$-height $n$ over the algebraic closure $\bbF$
of the residue field of $\fro$. Generalizing work of Lubin and
Tate, cf. \cite{LT}, V. G. Drinfeld showed in \cite{D} that the
functor of deformations of $\bbX$ is representable by an affine
formal scheme $\cM_0 = \Spf(R)$, where $R \simeq
\fronr[[u_1,\ldots,u_{n-1}]]$. Here, $\fronr$ is the completion of
the maximal unramified extension of $\fro$. Moreover, Drinfeld
introduced the notion of a level-$m$-structure, and proved that
the functor of deformations of $\bbX$ which are equipped with a
level-$m$-structure is representable by a formal scheme $\cM_m =
\Spf(R_m)$, where $R_m$ is a regular local ring which is a finite
flat $R$-module, and \'etale over $R$ after inverting a
uniformizer $\vpi$ of $\fro$. In this paper we are interested in
the geometrically connected components of the rigid-analytic space

$$M_m = \cM_m^{rig}$$

\medskip

associated to $\cM_m$, cf. \cite{dJ2}, sec. 7. The \'etale
cohomology groups of the spaces $M_m$ have been investigated in
the last two decades by various authors because of their
significance for the local Langlands correspondence. Let us cite
only \cite{Ca}, \cite{Bo}, and \cite{HT}. According to conjectures
by Carayol and Drinfeld, the inductive limit

$$H^{n-1} = \lim_{\stackrel{\lra}{m}} H^{n-1}(M_m
\times_\hFnr \bbC_\vpi,\Qlb)$$

\medskip

realizes simultaneously the Jacquet-Langlands and the Langlands
correspondence (cf. \cite{Ca} for a more precise statement). Here
$\hFnr$ is the field of fractions of $\fronr$, and $\bbC_\vpi$ is
the completion of an algebraic closure of $\hFnr$.\\

Whereas the spaces $M_m$ are defined purely locally, the analysis
of the inductive limit above is carried out in \cite{Bo} and
\cite{HT} by embedding the local situation into a global one. This
is done because it is very hard to understand the action of the
inertia group on $H^{n-1}$ (however, there are results for $m=1$
or $n=2$, cf. \cite{Y} resp. \cite{W}). By studying the geometry
of the spaces $M_m$ purely locally, it is possible to understand
the action of the pro-covering group $GL_n(\fro)$ and the action
of $\Aut(\bbX)$, thereby proving the assertion concerning the
Jacquet-Langlands correspondence (cf. \cite{St}, where the results
are unconditionally proved for the Euler-Poincar\'e characteristic
of the cohomology). Of course, one would like to understand
explicitly the cohomology in all degrees, and the most basic task
is hence to determine the representation on $H^0(M_m \times_\hFnr
\bbC_\vpi,\Qlb)$. This means to identify the group action on the
set of geometrically connected components of $M_m$. \\

We are going to explain how we do this. For $m \ge 0$ put $K_m =
Frac(R_m)$, $K = K_0$, and fix a separable closure $K^s$ of $K$
containing all $K_m$. Denote by $G_K = Gal(K^s/K)$ the absolute
Galois group of $K$. Let $X^{univ}$ be the universal deformation
of $\bbX$ over $R$, and let

$$T_{\vpi^m} = X^{univ}[\vpi^m](K^s)$$

\medskip

be the group of $K^s$-valued points of the $\vpi^m$-torsion
subgroup of $X^{univ}$. $T_{\vpi^m}$ is free of rank $n$ over
$\fro/(\vpi^m)$. Then we show, exactly as Raynaud in \cite{R},
that the induced action of $G_K$ on $\Lambda^n(T_{\vpi^m})$
factorizes through the canonical isomorphism

$$Gal(\hFnr_m/\hFnr) \lra (\fro/(\vpi^m))^\times \,,$$

\medskip

where $\hFnr_m \sub K^s$ is obtained from $\hFnr$ by adjoining all
$\vpi^m$-torsion points of a fixed formal $\fro$-module $LT$ of
height one over $\fronr$. As $Gal(K^s/K_m)$ acts trivially on
$T_{\vpi^m}$, and hence on $\Lambda^n(T_{\vpi^m})$, this implies
that $\hFnr_m$ is contained in $K_m$. Then one shows that
$R_m/\vpi_m$ is reduced, where $\vpi_m$ is a uniformizer of
$\hFnr_m$, and using a result of de Jong, cf. \cite{dJ2}, 7.3.5,
one obtains that $M_m$ is geometrically connected over $\hFnr_m$. \\

When thinking about this problem of geometrically connected
components we were inspired by de Jong's paper \cite{dJ1}. De Jong
also uses the crucial fact that the action on the determinant of
the Tate module is given by the cyclotomic character (he considers
only the case $F = \bbQ_p$), but his further reasoning is
different from ours as he is interested in a description of the
\'etale fundamental group of the corresponding period space,
cf. \cite{dJ1}, Prop. 7.4.\\

We conclude with a remark on an earlier approach to the problem of
geometrically connected components of these spaces. Let $\fronr_m$
be the ring of integers of the Lubin-Tate extension $\hFnr_m$.
From the inclusion $\fronr_m \sub R_m$ we get a morphism of formal
schemes

$$\cM_m \lra \cM_m^{(1)} := \Spf(\fro_{\hFnr_m}) \,.$$

\medskip

$\cM_m^{(1)}$ can be intrinsically defined as the deformation
space with level-$m$-structures of the reduction $LT_\bbF$ of
$LT$. Our original aim was to define {\it a priori} a functorial
map like this, and then to deduce that $M_m$ is geometrically
connected over $(\cM_m^{(1)})^{rig} = \Sp(\hFnr_m)$. A functorial
map as above can be thought of as associating to a deformation of
$\bbX$ which is equipped with a level-$m$-structure:

$$(X, (\vpi^{-m}\fro/\fro)^n \stackrel{\phi}{\lra} X[\vpi^m])$$

\medskip

its {\it determinant}:

$$(\Lambda^n(X), \Lambda^n(\vpi^{-m}\fro/\fro)^n
\stackrel{\Lambda^n\phi}{\lra} \Lambda^n(X)[\vpi^m]) \,.$$

\medskip

In an unpublished manuscript J. Lubin treats the problem of
defining the determinant of a one-dimensional formal module
together with a determinant map, cf. \cite{L}, and we, in our
earlier approach, worked along the same lines. First one has to
define $\Lambda^n(X)$, what one may do using, for example, Zink's
theory of displays (resp., in the equal characteristic case, the
'module des coordonn\'ees', cf. \cite{Bo}). Then, and more
difficult it seems, one has to compare the Tate modules. Using
Falting's results in \cite{F}, cf. p. 278, one may possibly do
this. The approach of this paper, however, seems
to be much more elementary.\\

{\it Acknowledgements.} I would like to thank Professor J. Lubin
for sending me his unpublished manuscript on determinants of
formal groups. Moreover, I thank L. Fargues for very helpful
discussions on the subject of this paper.\\

{\bf Notation.} In this paper, $F$ will be a non-Archimedean local
field, with ring of integers $\fro$, and $\vpi$ will be a
uniformizer of $F$. The number of elements of the residue field
will be denoted by $q$, and the residue field itself by $\bbF_q$.
We denote by $\bbF$ an algebraic closure of $\bbF_q$. $\hFnr$ is
the completion of the maximal unramified extension of $F$, and
$\fronr$ its ring of integers. $\bbC_\vpi$ denotes a the
completion of an algebraic closure of $\hFnr$. If $A$ is a local
ring we denote by $\frm_A$ its maximal ideal. The residue field of
a point $x$ on a scheme
will be denoted by $\kappa(x)$.\\

\section{Preliminaries}

In this section we recall without proof some facts about the
formal deformation schemes $\cM_m$, cf. \cite{D}, \cite{St}.\\

\begin{para}\label{deformations}
Let $\bbX$ be a one-dimensional formal group over $\bbF$ that is
equipped with an action of $\fro$, i.e. we assume given a
homomorphism $\fro \ra \End_\bbF(\bbX)$ such that the action of
$\fro$ on the tangent space is given by the reduction map $\fro
\ra \bbF_q \sub \bbF$. Such an object is  called a {\it formal
$\fro$-module} over $\bbF$. Moreover, we assume that $\bbX$ is of
$F$-height $n$, which means that the kernel of multiplication by
$\vpi$ is a finite group scheme of
rank $q^n$ over $\bbF$.\\

It is known that for each $n \in \bbZ_{>0}$ there exists a formal
$\fro$-module of $F$-height $n$ over $\bbF$, and that it is unique
up to isomorphism \cite{D}, Prop. 1.6, 1.7.\\

Let $\cC$ be the category of complete local noetherian
$\fronr$-algebras with residue field $\bbF$. A {\it deformation}
of $\bbX$ over an object $A$ of $\cC$ is a pair $(X,\iota)$,
consisting of a formal $\fro$-module $X$ over $A$ which is
equipped with an isomorphism $\iota: \bbX \ra X_{\bbF}$ of formal
$\fro$-modules over $\bbF$, where $X_{\bbF}$ denotes the reduction
of $X$ modulo the maximal ideal $\frm_A$ of $A$. Sometimes we will
omit $\iota$ from the notation.\\

Following Drinfeld \cite{D}, sec. 4B, we define a {\it structure
of level m} or {\it level-m-structure} on a deformation $X$ over
$A \in \cC$ ($m \ge 0$) as an $\fro$-module homomorphism

$$\phi: (\vpi^{-m}\fro/\fro)^n \lra \frm_A \,,$$

\medskip

such that, after having fixed a coordinate $T$ on the formal group
$X$, the power series $[\vpi]_X(T) \in A[[T]]$, which describes
the multiplication by $\vpi$ on $X$, is divisible by

$$ \prod_{a \in (\vpi^{-1}\fro/\fro)^n} (T-\phi(a)) \,.$$

\medskip

Here, $\frm_A$ is given the structure of an $\fro$-module via $X$.\\

Define the following set-valued functor $\cM_m$ on the category
$\cC$. For an object $A$ of $\cC$ put

$$\cM_m(A)=\{(X,\iota,\phi) \,|\, (X,\iota)
\mbox{ is a def. over $A$, $\phi$ is a level-$m$-structure on $X$}
\}/ \simeq \,,$$

\medskip

where $(X,\iota,\phi) \simeq (X',\iota',\phi')$ if and only if
there is an isomorphism $(X,\iota) \ra (X',\iota')$ of formal
$\fro$-modules over $A$, which is compatible with the level
structures. For $0 \le m' \le m$ one gets by restricting a
level-$m$-structure to

$$(\vpi^{-m'}\fro/\fro)^n \sub (\vpi^{-m}\fro/\fro)^n$$

\medskip

a level-$m'$-structure and hence a natural transformation

$$\cM_m \lra \cM_{m'} \,.$$

\medskip

Put $\fro_B = \End_\fro(\bbX)$. $\fro_B$ is the ring of integers
in a central division algebra over $F$ with Hasse invariant
$\frac{1}{n}$ (\cite{D}, Prop. 1.7). There is a natural action of
$GL_n(\fro/\vpi^m) \times \froxB$ from the right on the functor
$\cM_m$ given by

$$[X,\iota,\phi] \cdot (g,b) = [X,\iota \circ b,\phi \circ g]$$

\medskip

where $(g,b) \in GL_n(\fro/\vpi^m) \times \froxB$ and
$[X,\iota,\phi]$ denotes the equivalence class of $(X,\iota,\phi)$.\\
\end{para}

\begin{thm}\label{representability of def functors}
(i) The functor $\cM_m$ is representable by a regular local ring
$R_m$ of dimension $n$. Hence there is a universal formal
$\fro$-module $X^{univ}$ over $R := R_0$ which defines on the
maximal ideal $\frm_{R_m}$ of $R_m$ the structure of an
$\fro$-module. There is a universal level-m-structure

$$\phi^{univ}_m: (\vpi^{-m}\fro/\fro)^n \lra \frm_{R_m}$$

\medskip

such that, if $a_1,...,a_n$ is a basis of $(\vpi^{-m}\fro/\fro)^n$
over $\fro/(\vpi^m)$, then

$$\phi^{univ}_m(a_1), \ldots ,\phi^{univ}_m(a_n)$$

\medskip

is a regular system of parameters for $R_m$.

(ii) The ring homomorphism $R_m \ra R$ which corresponds to the
natural transformation $\cM_m \ra \cM_0$ makes $R_m$ a finite and
flat $R$-algebra. Moreover, $R_m[\frac{1}{\vpi}]$ is \'etale and
galois over $R[\frac{1}{\vpi}]$ with Galois group isomorphic to
$GL_n(\fro/\vpi^m)$.

(iii) $R$ is (non-canonically) isomorphic to
$\fronr[[u_1,\ldots,u_{n-1}]]$.

\end{thm}

{\it Proof.} (i) This result is \cite{D}, Prop. 4.3.\\

(ii) That $R_m$ is finite and flat over $R$ is again \cite{D},
Prop. 4.3. For the second statement we refer to \cite{St}, Thm. 2.1.2.\\

(iii) This is \cite{D}, Prop. 4.2. \hfill $\Box$ \\

{\it Remark.} The fact that $\fronr[[u_1,\ldots,u_{n-1}]]$
represents $\cM_0$ is due to Lubin and Tate (for $F=\bbQ_p$), cf.
\cite{LT}. For this reason $\cM_0$, the deformation space without
level structures, is sometines called the {\it Lubin-Tate moduli
space}, cf. \cite{HG}, \cite{Ch}.\\

By the preceding theorem, $R_m$ is a domain, and we put $K_m =
Frac(R_m)$. This is a Galois extension of $K := K_0$ with Galois
group canonically isomorphic to $GL_n(\fro/\vpi^m)$. The maps
$\cM_m \ra \cM_{m'}$ induce injections $K_{m'} \hra K_m$. Put

$$K_\infty = \cup_{m \ge 0} K_m \,,$$

\medskip

and fix a separable closure $K^s$ of $K$ containing
$K_\infty$. \\

\begin{para}\label{multiplication by vpi}
We conclude this section by recalling that one may choose the
parameters $\vpi = u_0, u_1,\ldots,u_{n-1}$ of $R$ and the
coordinate on $X^{univ}$ such that multiplication by $\vpi$ on
$X^{univ}$ is given by a power series $[\vpi]_{X^{univ}}(T) \in
R[[T]]$ with the property that

$$[\vpi]_{X^{univ}}(T) \equiv u_iT^{q^i} \mbox{ mod }
(u_0,\ldots,u_{i-1}) \,,  \mbox{ deg}(q^i + 1) \,,$$

\medskip

cf. \cite{HG}, Prop. 5.7. In particular, if $x \in \Spec(R)$ is a
point where $\vpi$ vanishes but $u_1$ is invertible in the residue
field $\kappa(x)$ of $x$, the multiplication of $\vpi$ on the
formal group $X^{univ} \hat{\otimes} \kappa(x)$ has as kernel a
group scheme of order $q$. Therefore the connected component of
the associated $\vpi$-divisible group $(X^{univ}[\vpi^\infty])
\otimes \kappa(x)$ over $\kappa(x)$ is a formal $\fro$-module of
height one.\\
\end{para}

\section{The Galois action on the determinant of the Tate module}

\begin{para}\label{Tate module}
Denote by $T$ the Tate module of the $\vpi$-divisible group
$X^{univ}[\vpi^\infty] \otimes K$:

$$T = \lim_{\stackrel{\longleftarrow}{m}} X^{univ}[\vpi^m](K^s) \,.$$

\medskip

The universal Drinfeld level-structures furnish an isomorphism of
$\fro$-modules

$$\fro^n \lra T \,,$$

\medskip

so that $\Lambda^n_\fro(T)$ is free of rank one over $\fro$.\\
\end{para}

\begin{para}\label{height one}
We recall some facts from Lubin-Tate theory. Fix a formal
$\fro$-module $LT$ of height one over $\fronr$. As the universal
deformation ring of height one formal $\fro$-modules is just
$\fronr$, all such formal $\fro$-modules are isomorphic. Denote by
$\hat{F}^{nr,s} \sub K^s$ the algebraic closure of $\hFnr$ in
$K^s$. It is a separable closure of $\hFnr$. Let $\frm$ be the
maximal ideal of the ring of integers in the completion of
$\hat{F}^{nr,s}$. $\frm$ is equipped via $LT$ with an
$\fro$-module structure, and the torsion points of $LT$ in $\frm$
are known to lie in $\hat{F}^{nr,s}$. Let $\hFnr_m \sub
\hat{F}^{nr,s}$ be the subfield generated over $\hFnr$ by the
$\vpi^m$-torsion points of $LT$ in $\frm$. As all $\fro$-modules
of height one over $\fronr$ are isomorphic over $\fronr$, this
field is independent of the choice of $LT$. There is a canonical
isomorphism

$$\chi_m: Gal(\hFnr_m/\hFnr) \lra (\fro/(\vpi^m))^\times \,,$$

\medskip

such that for any $\vpi^m$-torsion point $\alpha$ of $LT$ and
$\sigma \in Gal(\hFnr_m/\hFnr)$ one has

$$\sigma(\alpha) = [\chi_m(\sigma)]_{LT}(\alpha) \,.$$

\medskip

The field $\hFnr_\infty = \bigcup_m \hFnr_m$ is the maximal
abelian extension of $\hFnr$. The characters $\chi_m$ induce an
isomorphism

$$\chi: Gal(\hFnr_\infty/\hFnr) \lra \fro^\times \,,$$

\medskip

Finally, let $\tilde{\chi}: G_K \ra \fro^\times$ be the
composition of $G_K \ra Gal(\hFnr_\infty/\hFnr)$ with $\chi$.\\
\end{para}

\bigskip

The following theorem is essentially Raynaud's theorem on the
action of the Galois group on the determinant of the Tate module
of a $p$-divisible group, cf. \cite{R}, Thm. 4.2.1. We prove it
here again for the sake of completeness and because we need it in
the more general
context of formal $\fro$-modules.\\

\begin{thm}\label{galois on determinant}
The natural action of $G_K = Gal(K^s/K)$ on $\Lambda^n_\fro(T)$ is
given by the character $\tilde{\chi}$, i.e. for all $\sigma \in
G_K$ and $\lambda \in \Lambda^n_\fro(T)$ one has:

$$\sigma(\lambda) =
\tilde{\chi}(\sigma) \lambda \,.$$

\end{thm}

{\it Proof.} We assume $n \ge 2$, because for $n=1$ everything is
trivial. We follow the reasoning of the proof of \cite{R}, Th.
4.2.1. Let $x \in \Spec(R)$ be a prime ideal containing $\vpi$ but
not containing $u_1$, and denote by $\kappa(x)$ the residue field
at $x$. We already noticed in \ref{multiplication by vpi} that the
connected component of the $\vpi$-divisible group
$X^{univ}[\vpi^\infty] \otimes \kappa(x)$ is a formal
$\fro$-module of height one, and the \'etale part is hence of
height $n-1$. Let $R_x^{sh}$ be the strict henselization of $R$ at
$x$. Its residue field $\kappa(x)^s$ is a separable closure of
$\kappa(x)$. By \cite{D}, Prop. 1.7, all formal $\fro$-modules of
finite height over a separably closed field are isomorphic, so
that the formal module associated to the connected component

$$(X^{univ}[\vpi^\infty] \otimes \kappa(x)^s)^\circ$$

\medskip

is isomorphic to $LT \hat{\otimes}_\fronr \kappa(x)^s$. Therefore,
the formal module associated to

$$(X^{univ}[\vpi^\infty] \otimes R_x^{sh})^\circ$$

\medskip

is isomorphic to a deformation of $LT \hat{\otimes}_\fronr
\kappa(x)^s$, hence it is isomorphic to

$$LT \hat{\otimes}_\fronr R_x^{sh} \,,$$

\medskip

because there is up to isomorphism only one deformation, cf.
\cite{Ha}, Thm. 22.4.16. Let $R_x^h \sub R_x^{sh}$ be the
henselization of $R$ at $x$, put $K_x^h = Frac(R_x^h)$, and fix a
separable closure $K_x^{h,s}$ of $K_x^h$ together with an
embedding of $K^s$ into this separable closure. We get an
embedding of Galois groups $Gal(K_x^{h,s}/K_x^h) \hra G_K$, and it
follows from what we have said above that the restriction of the
$G_K$-action on

$$\Pi := \Lambda^n_\fro(T) \otimes \tilde{\chi}^{-1}$$

\medskip

to $Gal(K_x^{h,s}/K_x^h)$ factors through the absolute Galois
group of $\kappa(x)$. That means, the representation of $G_K$ on
$\Pi$ is unramified at $x$. As the $\vpi$-divisible group
$X^{univ}[\vpi^\infty]$ is \'etale at every point of $\Spec(R)$
where $\vpi$ is invertible, $\Pi$ is unramified at all points of
$U = \Spec(R) - V((\vpi,u_1))$, whose complement is of codimension
$2$. Hence the representation of $G_K$ on $\Pi$ extends to the
\'etale fundamental group $\pi_1(U,\Spec(K^s))$ of $U$. But by the
Zariski-Nagata purity theorem, cf. \cite{SGA2}, exp. X, Thm. 3.4,
\cite{SGA1}, exp. X, Cor. 3.3, $\pi_1(U,\Spec(K^s)) =
\pi_1(\Spec(R),\Spec(K^s))$, and the latter group is trivial, as
$R$ is strictly henselian. This means that $G_K$ acts trivially on
$\Pi$ and therefore by $\tilde{\chi}$ on $\Lambda^n_\fro(T)$.
\hfill $\Box$

\begin{cor}\label{galois action finite level}
(i) The action of $G_K$ on $\Lambda^n_\fro(X^{univ}[\vpi^m](K^s))$
is given by the character $\tilde{\chi}_m$, which is the
composition of $\tilde{\chi}$ with the canonical map $\fro^\times
\ra (\fro/\vpi^m)^\times$.\\

(ii) The ring $R_m = {\mathcal O}(\cM_m)$ contains the ring of
integers $\fronr_m$ of the Lubin-Tate extension $\hFnr_m$.
\end{cor}

{\it Proof.} (i) This assertion follows immediately from the
$G_K$-equivariant isomorphism

$$\Lambda^n_\fro(X^{univ}[\vpi^m](K^s)) \simeq
\Lambda^n_\fro(T)/\vpi^m\Lambda^n_\fro(T) \,.$$

\medskip

(ii) The subgroup $Gal(K^s/K_m) \sub G_K$ acts trivially on
$X^{univ}[\vpi^m](K^s)$, and hence trivially on
$\Lambda^n(X^{univ}[\vpi^m](K^s))$. By (i) the action on the
latter module is given by $\tilde{\chi}_m$. Therefore,
$Gal(K^s/K_m)$ acts trivially on $\hFnr_m$, and hence $\hFnr_m$ is
contained in $K_m$. Because $R_m$ is integrally closed, it
contains $\fronr_m$. \hfill $\Box$

\medskip

By the second assertion, we will view from now on $R_m$ as an
$\fronr_m$-algebra, and $R_m[\frac{1}{\vpi}]$ as an
$\hFnr_m$-algebra.\\

\bigskip

\section{Geometrically connected components and the group action on $\pi_0$}

\begin{para}
Fix an integer $m \ge 1$. In this section we show first that
$R_m[\frac{1}{\vpi}]$ is geometrically integral over $\hFnr_m$.
Then we determine the structure of $\pi_0(M_m \times_\hFnr
\bbC_\vpi)$ together with the action of
$GL_n(\fro) \times \froxB \times G_\hFnr$.\\

Denote as in Cor. \ref{galois action finite level} by $\fronr_m$
the ring of integers in $\hFnr_m$, and let $\vpi_m$ be a
uniformizer in $\fronr_m$. Furthermore, we fix a set $\cR \sub
(\vpi^{-m}\fro/\fro)^n$ of representatives of the orbits of the
action of $(\fro/\vpi^m)^\times$ on

$$(\vpi^{-m}\fro/\fro)^n - (\vpi^{-(m-1)}\fro/\fro)^n \,.$$

\medskip

Finally, we abbreviate the universal level-$m$-structure
$\phi^{univ}_m$, cf. \ref{representability of def functors}, by
$\phi_m$.\\
\end{para}

\begin{prop}\label{integralness}
(i) For $\al \in \cR$ the element $\phi_m(\al) \in \frm_{R_m}$ of
the regular local ring $R_m$ is irreducible. Moreover, the
elements $\phi_m(\al)$ et $\phi_m(\beta)$ are not associated if
$\al \neq \beta$, $\al, \beta \in \cR$.\\
(ii) Up to a unit in $R_m$, the element $\vpi_m$ is equal to the product
$\prod_{\al \in \cR} \phi_m(\al)$.\\
(iii) The ring $R_m/\vpi_mR_m$ is reduced.\\
(iv) The ring $R_m[\frac{1}{\vpi}]$ is geometrically integral over
$\hFnr_m$.
\end{prop}

{\it Proof.} (i) An element $\al \in\cR$ can be completed to a
basis $\al = \al_1, \ldots, \al_n$ of $(\vpi^{-m}\fro/\fro)^n$. By
\ref{representability of def functors} we know that
$\phi_m(\al_1),\ldots,\phi_m(\al_n)$ is a regular set of
parameters of $R_m$. Hence, the element $\phi_m(\al) =
\phi_m(\al_1)$ is in particular irreducible. Let $\al, \beta$ be
two different elements of $\cR$. For a point $x \in \Spec(R_m)$
consider the induced homomorphism

$$\phi_{m,x}:(\vpi^{-m}\fro/\fro)^n \stackrel{\phi_m}{\lra}
X^{univ}[\vpi^m](R_m) \lra X^{univ}[\vpi^m](\kappa(x)) \,.$$

\medskip

If $h$ is the height of the connected part of the $\vpi$-divisible
group $X^{univ}[\vpi^\infty] \otimes \kappa(x)$, the kernel of
$\phi_{m,x}$ is a direct summand of $(\vpi^{-m}\fro/\fro)^n$ of
rank $h$ over $\fro/(\vpi^m)$. And conversely, if the rank of the
kernel of $\phi_{m,x}$ is equal to $h$, then the height of the
connected component of $X^{univ}[\vpi^\infty] \otimes \kappa(x)$
is $h$. Hence, if $\al$ and $\beta$ are two different elements of
$\cR$, the kernel of $\phi_{m,x}$ is at least of rank two if it
contains $\al$ and $\beta$, and in this case the height of the
connected component of $X^{univ}[\vpi^\infty] \otimes \kappa(x)$
is at least two. Now suppose that $\phi_m(\al)$ and
$\phi_m(\beta)$ are associated prime elements. Let $x_0 \in
\Spec(R_0)$ correspond to the prime ideal $\vpi R_0$. Then, as we
already stated in \ref{multiplication by vpi}, the height of the
connected component of $X^{univ}[\vpi^\infty] \otimes \kappa(x_0)$
is one. As $R_m$ is finite and flat over $R_0$, there is a prime
ideal $x_m$ of $R_m$ lying over $x_0$. Then the kernel of
$\phi_{m,x_m}$ is a direct summand of rank one, hence generated by
one element, $\gamma$ say. Let $g \in GL_n(\fro/\vpi^m)$ be an
element with $g(\al) = \gamma$. Then $\al$ lies in the kernel of
$\phi_{m,y_m}$, where $y_m = g(x_m)$, so $\phi_m(\al)$ is in the
prime ideal corresponding to $y_m$. But if $\phi(\al)$ and
$\phi(\beta)$ are associated, i.e. $\phi(\beta) = u \phi(\al)$
with a unit $u \in R_m$, we also have $\phi_{m,y_m}(\beta)=0$, and
so the kernel of $\phi_{m,y_m}$ is at least of rank two, which cannot be,
because $\phi_{m,y_m} = \phi_{m,x_m} \circ g$.\\

(ii) Let $\al \in \cR$. Then $\frp = \phi(\al)R_m$ is a prime
ideal, and because the kernel of $\phi_{m,\frp}$ contains $\al$,
the connected component of $X^{univ}[\vpi^\infty] \otimes
\kappa(\frp)$ is at least of height one. Hence $\frp$ contains
$\vpi$, and therefore $\vpi_m$ too. Because the elements
$\phi_m(\al)$, $\al \in \cR$, are pairwise not associated, and
because $R_m$ is a UFD, we see that

$$\vpi_m = f \cdot \prod_{\al \in \cR} \phi(\al) \,,$$

\medskip

with some element $f \in R_m$. Now let us consider the unramified
extension $F'$ of $F$ of degree $n$ inside $\hFnr$. Denote by
$\fro'$ the ring of integers of $F'$ and fix a formal
$\fro'$-module $X$ of $F'$-height one. Let $\iota: \bbX \ra X
\times_{\fro'} \bbF$ be an isomorphism, such that the pair $(X
\times_{\fro'} \fronr, \iota)$ is a deformation of $\bbX$ (as a
formal $\fro$-module of height $n$). This pair corresponds to a
point $x$ in $\Spec(R_0)$. We lift this point to a point $y$ in
$\Spec(R_m)$. The residue field $\kappa(y)$ at $y$ is then an
extension of $F'$ generated by the $\vpi^m$-torsion points of $X$.
Denote by $v$ the valuation on this extension which is normalized
by $v(\vpi) = 1$. Then the valuation of a $\vpi^m$-torsion point
of $X$, which is not annihilated by $\vpi^{m-1}$, is equal to
$\frac{1}{(q^n-1)q^{n(m-1)}}$. Mapping $\vpi_m$ into $\kappa(y)$
and using the equation above we calculate

$$\begin{array}{rl}
\frac{1}{(q-1)q^{m-1}} & = v(\vpi_m) \ge \sum_{\al \in
\cR}v(\phi_{m,y}(\al))
= \sum_{\al \in\cR}\frac{1}{(q^n-1)q^{n(m-1)}} \\
 & \\
 & = \frac{(q^n-1)q^{n(m-1)}}{(q-1)q^{m-1}}\frac{1}{(q^n-1)q^{n(m-1)}}
 = \frac{1}{(q-1)q^{m-1}} \,.
\end{array}$$\\

\medskip

And this shows that $f$ is necessarily a unit in $R_m$.\\

(iii) This is an immediate consequence of (i) and (ii).\\

(iv) By \cite{EGA}, Cor. 18.9.8, it suffices to show that the
fibres of

$$\Spec(R_m) \lra \Spec(\fronr_m)$$

\medskip

are geometrically reduced. By (iii), this is the case for the
fibre over the closed point. Let us now consider the generic
fibre. Let $E$ be a field containing $\hFnr_m$. Because
$R_m[\frac{1}{\vpi}]$ is \'etale over $R_0[\frac{1}{\vpi}]$, the
ring extension

$$R_0[\frac{1}{\vpi}] \otimes_\hFnr E \ra R_m[\frac{1}{\vpi}] \otimes_\hFnr E
= \prod_{\sigma \in Gal(\hFnr_m/\hFnr)}R_m[\frac{1}{\vpi}]
\otimes_{\hFnr_m} E$$

\medskip

is \'etale too. Because $R_0[\frac{1}{\vpi}] \otimes_\hFnr E$ maps
injectively into $E[[u_1, \ldots, u_{n-1}]]$, this ring is
reduced. By general results on \'etale extensions,

$$R_m[\frac{1}{\vpi}] \otimes_{\hFnr_m} E \,,$$

\medskip

which is \'etale over $R_0[\frac{1}{\vpi}] \otimes_\hFnr E$, is
reduced too. Therefore $R_m[\frac{1}{\vpi}] \otimes_{\hFnr_m} E$
is reduced. \hfill $\Box$

\bigskip

\begin{thm}\label{connected components}
(i) Let $E$ be a finite separable extension of $\hFnr$, which
contains the Lubin-Tate extension $\hFnr_m$ of $\hFnr$. Then the
rigid-analytic space

$$M_m \times_\hFnr \Sp(E) = \cM_m^{rig} \times_\hFnr \Sp(E)$$

\medskip

over $E$ has $(q-1)q^{m-1}$ connected components. These are the
fibres of the morphism

$$M_m \times_\hFnr \Sp(E) \lra \Sp(\hFnr_m) \times_\hFnr \Sp(E)
= \Sp(\hFnr_m \otimes_\hFnr E) \,.$$

\medskip
\end{thm}

{\it Proof.} For the construction of the rigid-analytic space
associated to a formal scheme we refer to \cite{dJ2}, sec. 7. The
first assertion clearly follows from the second. Because

$$M_m \times_\hFnr \Sp(E) =  \coprod_{\sigma \in Gal(\hFnr_m/\hFnr)}
M_m \times_{\hFnr_m} \Sp(E) \,,$$

\medskip

we only need to show that $M_m \otimes_{\hFnr_m} \Sp(E)$ is
connected. By \cite{dJ2}, 7.2.4 (g), we have

$$M_m \otimes_{\hFnr_m} \Sp(E) = \Spf\left(R_m \hat{\otimes}_{\fronr_m}
\fro_E\right)^{rig} \,,$$

\medskip

where $\fro_E$ is the ring of integers in $E$. By \cite{dJ2},
7.3.5, this space is connected if the ring

$$R_m \hat{\otimes}_{\fronr_m} \fro_E = R_m \otimes_{\fronr_m}
\fro_E$$

\medskip

(this equality holds because $E/\hFnr_m$ is finite) is integrally
closed. This ring is contained in $R_m \otimes_{\fronr_m} E =
R_m[\frac{1}{\vpi}] \otimes_{\hFnr_m} E$ which is integral, by the
preceding proposition \ref{integralness}. By \cite{Bou}, V, \S
1.7, Cor. to Prop. 19, $R_m \otimes_{\fronr_m} \fro_E$ is
integrally closed. \hfill $\Box$.

\bigskip

Let $\bbC_\vpi$ be a completion of an algebraic closure
$\bar{F}^{nr}$ of $\hFnr$. Denote by $G_\hFnr =
Gal(\bar{F}^{nr}/\hFnr)$ the absolute Galois group of $\hFnr$. By
continuity it acts on $\bbC_\vpi$. In the following we use the
isomorphism $\hFnr_m \otimes_\hFnr \bbC_\vpi \simeq \prod_{\sigma
\in Gal(\hFnr_m/\hFnr)} \bbC_\vpi$ given by

$$\lambda \otimes \mu \mapsto
(\sigma^{-1}(\lambda)\mu)_{\sigma \in Gal(\hFnr_m/\hFnr)} \,.$$

\medskip

This isomorphism is used to identify the connected components of
$\Sp(\hFnr_m \otimes_\hFnr \bbC_\vpi)$ with the connected
components of $\coprod_{\sigma \in Gal(\hFnr_m/\hFnr)}
\Sp(\bbC_\vpi)$ which we identify with its indexing set
$Gal(\hFnr_m/\hFnr)$. The latter group gets identified via the
character $\chi_m$ with $(\fro/\vpi^m)^\times$. Note that if we
let $\tau \in G_\hFnr$ act on the second factor of $\hFnr_m
\otimes_\hFnr \bbC_\vpi$, we have, via the isomorphism above,

$$((\sigma^{-1}(\lambda)\tau(\mu))_\sigma =
\tau((\sigma \circ \tau)^{-1}(\lambda)\mu)_\sigma) =
\tau((\sigma^{-1}(\lambda)\mu)_{\sigma \circ \tau^{-1}}) \,.$$

\medskip

Therefore, on the indexing set $Gal(\hFnr_m/\hFnr)$, $\tau$ acts
by multiplication by $(\tau|_{\hFnr_m})^{-1}$, and consequently on
$(\fro/\vpi^m)^\times$ by $\chi_m(\tau|_{\hFnr_m})^{-1}$.

\bigskip

\begin{thm}
(i) The morphism $M_m \ra \Sp(\hFnr_m)$ induces a bijection

$$\pi_0(M_m \times_\hFnr \Sp(\bbC_\vpi))
\stackrel{\simeq}{\lra} \pi_0(\Sp(\hFnr_m \times_\hFnr \bbC_\vpi))
\,,$$

\medskip

and the set on the right is identified with
$(\fro/\vpi^m)^\times$, as explained above. The resulting
bijection

$$\pi_0(M_m \times_\hFnr \Sp(\bbC_\vpi))
\stackrel{\simeq}{\lra} (\fro/\vpi^m)^\times$$

\medskip

is $GL_n(\fro) \times \froxB \times G_\hFnr$-equivariant if we let
$GL_n(\fro) \times \froxB \times G_\hFnr$ act on
$(\fro/\vpi^m)^\times$ by

$$(g,b,\tau) \mapsto \det(g) Nrd(b)^{-1}
\chi(\tau|_{\hFnr_\infty})^{-1} \,\, \mod \, (1+\vpi^m\fro) \,.$$

\medskip

Here $Nrd: \froxB \ra \fro^\times$ denotes the reduced norm. \\

(ii) In particular, the zero'th $l$-adic \'etale cohomology group
decomposes as follows:

$$H^0(M_m \times_\hFnr \Sp(\bbC_\vpi),\Qlb)
\, \simeq \, \bigoplus_\omega (\omega \circ \det) \otimes (\omega
\circ Nrd)^{-1} \otimes (\omega \circ rec_\hFnr) \,,$$

\medskip

where $\omega: (\fro/\vpi^m)^\times \ra \Qlb^\times$ runs through
all $\Qlb$-valued characters of $(\fro/\vpi^m)^\times$, and
$rec_\hFnr$ is the reciprocity map from local class field theory
(normalized such that an arithmetic Frobenius is mapped to a
uniformizer).\\
\end{thm}

{\it Proof.} (i) As $R_m$ is normal, by \cite{dJ2}, 7.3.5., the
rigid space $M_m = \Spf(R_m)^{rig}$ is connected. By \cite{Co},
Cor. 3.2.3, there is a finite separable extension $E$ of $\hFnr$
such that the canonical map

$$M_m \times_\hFnr \Sp(\bbC_\vpi) \lra M_m \times_\hFnr
\Sp(E)$$

\medskip induces a bijection of the connected components:

$$\pi_0(M_m \times_\hFnr \Sp(\bbC_\vpi)) \stackrel{\simeq}{\lra}
\pi_0(M_m \times_\hFnr \Sp(E)) \,.$$

\medskip

By theorem \ref{connected components}, we can take here $E =
\hFnr_m$, and get:

$$\begin{array}{rl}
\pi_0(M_m \times_\hFnr \Sp(\bbC_\vpi)) & = \pi_0(M_m \times_\hFnr
\Sp(\hFnr_m)) = \pi_0(\Sp(\hFnr_m
\otimes_\hFnr \hFnr_m))  \\
 & \\
 & =  \pi_0(\Sp(\prod_{\sigma \in Gal(\hFnr_m/\hFnr)} \hFnr_m)) \\
 & \\
 & = Gal(\hFnr_m/\hFnr) \stackrel{\chi}{\lra} (\fro/\vpi^m)^\times \,.

\end{array}$$

\bigskip

Here we used the isomorphism $\hFnr_m \otimes_\hFnr \hFnr_m \simeq
\prod_{\sigma \in Gal(\hFnr_m/\hFnr)} \hFnr_m$ given by

$$\lambda \otimes \mu \mapsto
(\sigma^{-1}(\lambda)\mu)_{\sigma \in Gal(\hFnr_m/\hFnr)} \,.$$

\medskip

The action of the Galois group $G_K$ on $R_m$ factors through
$GL_n(\fro)$, so let $g \in GL_n(\fro)$ come from some $\rho \in
G_K$. Then, for a $\vpi^m$-torsion point $a$ of $LT$ in $\hFnr_m$
we have

$$\rho(a) = [\tilde{\chi}(\rho)]_{LT}(a) \,.$$

\medskip

But $\tilde{\chi}(\rho)$ is the element by which $\rho$ acts on
the Tate module $\Lambda^n(T)$, by theorem \ref{galois on
determinant}. On the other hand, $g$ acts as $\det(g)$ on
$\Lambda^n(T)$. Therefore, $g$ acts on the torsion point $a \in
\fronr_m \sub R_m$ by $[\det(g)]_{LT}(a)$. Via the above
identification of $\pi_0(M_m \times_\hFnr \Sp(\bbC_\vpi))$ with
$(\fro/\vpi^m)^\times$, we get that $GL_n(\fro)$ acts on
$(\fro/\vpi^m)^\times$ via the determinant ($\mod
(1+\vpi^m\fro)$). The action of $G_\hFnr$ we already computed
above. Finally, we consider the action of $\froxB$. Let $\fro'
\sub \fro_B$ be the maximal unramified extension of $\fro$ in
$\fro_B$. Because the reduced norm on $\froxB$, when restricted to
$(\fro')^\times$, maps $(\fro')^\times$ surjectively onto
$\fro^\times$, it suffices to calculate the action of
$(\fro')^\times \sub \froxB$. Let $X$ be a formal $\fro'$-module
of height one over $\fro'$. Then $X \hat{\otimes}_{\fro'} \fronr$,
when equipped with an isomorphism of its special fibre with
$\bbX$, corresponds to a point $x$ in $\Spec(R_0)$ which is fixed
by the action of $(\fro')^\times$. Let $y \in \Spec(R_m)$ be a
point over $x$. Then $\fro'_m := R_m/y$ is the ring of integers in
$\hFnr.F'_m$, where $F'_m$ is the Lubin-Tate extension of $F' =
\fro'[\frac{1}{\vpi}]$ generated by the $\vpi^m$-torsion points of
$X$. By local class field theory, the action of an element $b \in
(\fro')^\times$ on the subfield $\hFnr_m \cap F'_m \sub F'_m$ is
given by the norm of $N_{F'/F}(b) = Nrd(b)$ (on torsion points of
$LT$ inside $\hFnr_m$). But the action of $b \in \froxB$ on a
triple $(X,\iota,\phi)$ is defined by $(X,\iota \circ b, \phi)$,
and the latter object is equivalent to $(X,\iota, b^{-1} \circ
\phi)$, because $X$ has multiplication by $\fro'$. Hence the
action of $b$ (via $R_m$ and specialisation) on the
$\vpi^m$-torsion points of $X$ is given by $[b^{-1}]_X$.
Therefore, it acts on the $\vpi^m$-torsion points of $LT$ by
$[Nrd(b)^{-1}]_{LT}$. Via the above identification of the
connected components of $M_m \times_\hFnr \Sp(\bbC_\vpi)$ with
$(\fro/\vpi^m)^\times$, $b$ acts on the latter set by
multiplication with $Nrd(b)^{-1}$. This proves that $\froxB$ acts
on $(\fro/\vpi^m)^\times$ by
$Nrd^{-1}$.\\

(ii) It is well known, that the local reciprocity map (normalized
such that an arithmetic Frobenius is mapped to a uniformizer)
fullfills

$$rec_\hFnr(\tau) = \chi(\tau)^{-1} \,.$$

\medskip

From this and the first part of the theorem the assertion follows
immediately. \hfill $\Box$.

\bigskip

We can finally also pass to the limit over all $m$ and get a
natural $GL_n(\fro) \times \froxB \times G_\hFnr$-equivariant map

$$\lim_{\stackrel{\longleftarrow}{m}} \pi_0(M_m \times_\hFnr \Sp(\bbC_\vpi))
\stackrel{\simeq}{\lra} \fro^\times \,,$$

\medskip

where on the right side the group $GL_n(\fro) \times \froxB \times
G_\hFnr$ acts by

$$(g,b,\tau) \mapsto \det(g) Nrd(b)^{-1}
rec_\hFnr(\tau) \,.$$

\medskip

And hence:

$$\lim_{\stackrel{\lra}{m}} H^0(M_m \times_\hFnr \Sp(\bbC_\vpi),\Qlb)
\, \simeq \, \bigoplus_\omega (\omega \circ \det) \otimes (\omega
\circ Nrd)^{-1} \otimes (\omega \circ rec_\hFnr) \,,$$

\medskip

where $\omega$ runs through all continuous characters $\fro^\times
\ra \Qlb^\times$ (with the discrete topology on $\Qlb^\times$).

\end{document}